\documentclass[11pt]{amsart}
\usepackage{amsmath}
\usepackage{amsthm}
\usepackage{amssymb}
\usepackage{amsfonts}
\usepackage{graphicx}
\usepackage{inputenc}
\usepackage{upref}
\usepackage{bm}
\usepackage{calc}
\usepackage{color}
\usepackage{ragged2e}

  \theoremstyle{plain}

\newtheorem{theorem}{Theorem}[section]

\newtheorem{definition}[theorem]{Definition}

\newtheorem*{theorem*}{Theorem}
\newtheorem*{definition*}{Definition}

\theoremstyle{definition}

\let\i\relax
\let\b\relax
\newcommand{\i}[1]{\textit{#1}}
\newcommand{\b}[1]{\textbf{#1}}

\let\int\relax
\newcommand{\int}[1]{\mathring{#1}}
\newcommand{\Z}{\mathbb{Z}}

\newcommand{\into}{\hookrightarrow}

\newcommand{\std}{\text{std}}
\newcommand{\pt}{\text{pt}}
\let\~\relax
\newcommand{\~}{\sim}
\usepackage[margin=1.5in]{geometry}

\let\int\relax
\newcommand{\int}{\mathring}

\title[Concordances from the standard surface in $S^2\times S^2$]{Concordances from the standard surface in $S^2\times S^2$}
\author[Maggie Miller]{Maggie Miller}
\thanks{The author is a fellow in the National Science Foundation Graduate Research Fellowship program, under Grant No. DGE-1656466.}
\address{Department of Mathematics, Princeton University, Princeton, NJ 08544}
\email{maggiem@math.princeton.edu}

\begin{document}
\maketitle
\begin{abstract}
In this note, we combine the recent 4-dimensional light bulb theorem of David Gabai and a recent construction of concordances for knots in $S^2\times S^1$ due to Eylem Zeliha Yildiz to construct a concordance between the standard surface of genus $g$ in $S^2\times S^2$ and any homologous surface.
\end{abstract}
\setcounter{tocdepth}{2}
\setcounter{equation}{0}
\section{Introduction}



Fix $x_0\in S^2$. We say the \emph{standard surface of genus-$g$} in $S^2\times S^2$ is a copy of $R^{\std}$ with $g$ standard handles attached,  where $R^{\std}$ is the sphere $x_0\times S^2$. This notation agrees with that of Gabai's paper on the $4$-dimensional light bulb theorem \cite{gabai}.

In this paper, we prove the following theorem.

\begin{theorem*}
Let $R$ be a connected orientable genus-$g$ surface embedded in $S^2\times S^2$ with $[R]=[R^{\std}]\in H_2(S^2\times S^2;\Z)$. Then $R$ is smoothly concordant to the standard genus-$g$ surface.
\end{theorem*}

We recall the standard definition of smooth concordance.
\begin{definition*}\rm
Let $X,Y$ be smooth submanifolds of a manifold $W$. We say $X$ and $Y$ are {\emph{concordant}} if there exists a smooth embedding $f:X\times I\into W\times I$ so that $f(X\times 0)=X\subset W\times 0=W$ and $f(X\times 1)=Y\subset W\times 1=W$.
\end{definition*}

We should note that the result follows as a consequence of Sunukjian's \cite{sunukjian} more general theorem on concordance in $4$-manifolds.
\newtheorem*{sunukjian*}{Theorem 6.1 of \cite{sunukjian}}
\begin{sunukjian*}
Any two homologous surfaces of the same genus within a simply connected $4$-manifold are smoothly concordant.
\end{sunukjian*}
Sunukjian's result uses the Thom construction and some surgery theory. We give this document as a simple construction of concordances between surfaces in $S^2\times S^2$, with the hope that this description might aid in the study of concordance of surfaces/spheres in arbitrary $4$-manifolds (perhaps without $2$-torsion in fundamental group; see the second half of \cite{gabai}).

Our proof follows from the following two results due respectively to Yildiz \cite{yildiz} and Gabai \cite{gabai}.

\newtheorem*{yildiz*}{Theorem 2 of \cite{yildiz}}
\newtheorem*{gabai*}{Theorem 9.8 of \cite{gabai}}
\newtheorem*{lightbulb*}{3-dimensional light bulb theorem}
\begin{yildiz*}\label{yildiz}
Let $K\subset S^2\times S^1$ be a knot transverse to $S^2\times \pt$ with algebraic intersection number $\langle K\cap S^2\times\pt\rangle=1$. Then $K$ is smoothly concordant to $\pt\times S^1$.
\end{yildiz*}

\begin{gabai*}\label{gabai}
Let $R$ be a connected embedded genus-$g$ surface in $S^2\times S^2$ such that $|R\cap (S^2\times y_0)|=1$ and $[R]=[R^{\std}]\in H_2(S^2\times S^2;\Z)$. Then $R$ isotopic to the standard genus-$g$ surface via an ambient isotopy that fixes $S^2\times y_0$ pointwise. 
\end{gabai*}

Gabai's theorem is a 4-dimensional analog of a classical result in dimension three.

\begin{lightbulb*}
Let $K\subset S^2\times S^1$ be a knot transverse to $S^2\times \pt$ with $| K\cap S^2\times\pt|=1$. Then $K$ is smoothly isotopic to $\pt\times S^1$ or $-(\pt\times S^1)$ (depending on the orientation of $K$).
\end{lightbulb*}

\subsection*{Acknowledgements}
The author would like to thank her advisor, David Gabai, for helpful suggestions and conversations (in particular about surfaces in $4$-dimensional space, but also more generally). Thanks also to Daniel Ruberman for pointing the author toward Sunukjian's theorem on concordance.

\begin{figure}
\includegraphics[width=\textwidth]{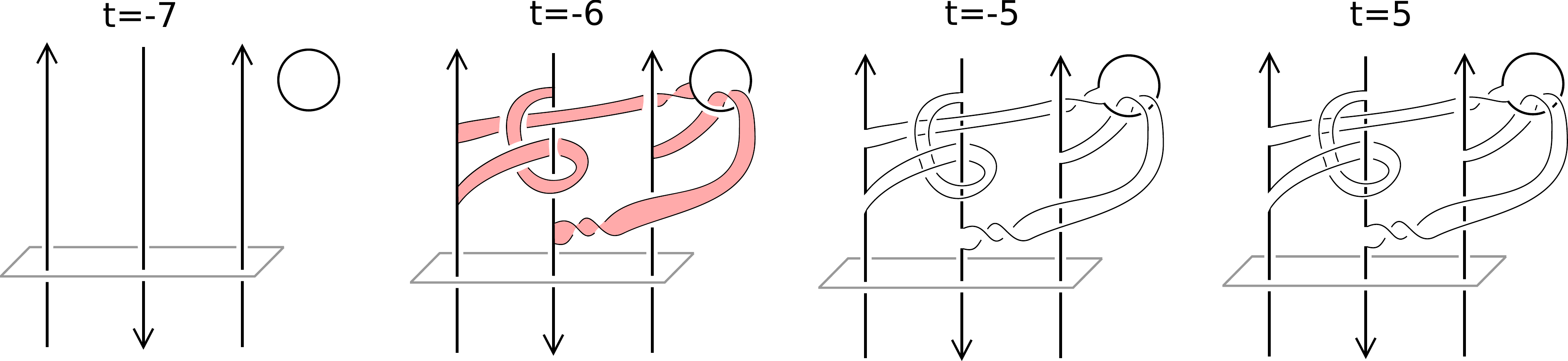}
\caption{A schematic of $R\cap \left(S^2\times S^1\times [-7,5]\right)$. For $-7\le t<-6$, $R\cap \left(S^2\times S^1\times t\right)$ consists of $2n+p+1$ curves. At time $t=-6$, bands appear. For $-6<t<6$, $R\cap\left( S^2\times S^1\times t\right)$ is a connected curve. 
The illustrations here of $t=-5,t=5$ are the same image, as $R\cap \left(S^2\times S^1\times (-6,6)\right)$ is taken to be a product.}
\label{fig:beforehandles}
\end{figure}

\section{Proof}

%
%


Yildiz's proof of concordance in dimension three is constructive. We will use this same construction in a neighborhood of a $3$-dimensional cross-section of $S^2\times S^2$, making use of Gabai's view of $S^2\times S^2$. All corners should be smoothed.

We position $R$ in $S^2\times S^2$ following the setup of Gabai \cite{gabai}.

\subsection*{Step 1. Fix a standard height function on $\mathbf{S^2\times S^2}$.}\leavevmode

This is (similar to) the normal form of \cite{suzuki}, as seen in \cite{gabai}.

Fix $y_0\in S^2$. Perturb $R$ so that $R$ is transverse to $S^2\times y_0$. Note the signed intersection $\langle R\cap\left( S^2\times y_0\right)\rangle=1$ and $R\cap (S^2\times y_0)$ contains $2n+1$ points.
View $S^2\times S^2$ as $S^2\times S^1\times[-\infty,\infty]/\~$, where $\~$  crushes each $S^1\times\pm\infty$ to a point. More explicitly, $(x,z,\pm\infty)\sim(x,z',\pm\infty)$ for each $x\in S^2$ and $z,z'\in S^1$.

 Here, $y_0\in S^2$ corresponds to the point $(z_0,0)\in S^1\times[\infty,\infty]/\~$ $\cong S^2$. We refer to the cross-section $S^2\times S^1\times t$ of $S^2\times S^2$ as being at height $t$.

Fix a small $\epsilon$ and isotope $R$ so that $R\cap (S^2\times y_0)$ is contained in a ball of radius $\epsilon$ about $0\times y_0$.
 Take $U$ to be the neighborhood of $S^2\times y_0$ given by $U=S^2\times \left(S^1\times\left([-\infty,-10)\cup(10,\infty]\right)\cup (-\epsilon,\epsilon)\times[-10,10]\right)$. We may assume $R\cap U$ is standard; i.e. $R\cap U=\left(\cup_{n+1} R^{\std}\cup_n-R^{\std}\right)\cap U=\cup\{2n+1$ disks$\}$.

While fixing $U$, $R$ can be isotoped to be transverse to each $S^2\times S^1\times t$ (for $t\neq-9,0,9$). At times $t=-9,0,9$, there are critical points of index $0,1,2$, respectively. We perturb $R$ so that these critical points are flattened into $2$-dimensional regions. Thus at time $t=-9$, $p$ disks (corresponding to local minima) appear. At time $t=0$, $4n+p+q+2g$ ribbon bands (corresponding to saddle points) appear. At $t=9$, $q$ disks (corresponding to local maxima) appear. 
Since $R$ is connected, there is some subset of the ribbon bands so that the result of resolving $R\cap \left(S^2\times S^1\times -8\right)$ along only these bands results in a connected curve. Push this subset of bands to height $t=-6$ and the rest to height $t=6$.

%
%

The current situation is partially illustrated in Figure \ref{fig:beforehandles}.

In particular, for $-6<t<6$, $R\cap(S^2\times S^1\times t)$ is a connected curve. Within the $3$-dimensional cross-section $S^2\times S^1\times t\cong S^2\times S^1$, this curve intersects $S^2\times (z_0,t)$ algebraically once.

\subsection*{Step 2. Attach handles in the region $\mathbf{-4\le t\le 4}$.}\leavevmode

\begin{figure}
\includegraphics[width=40mm]{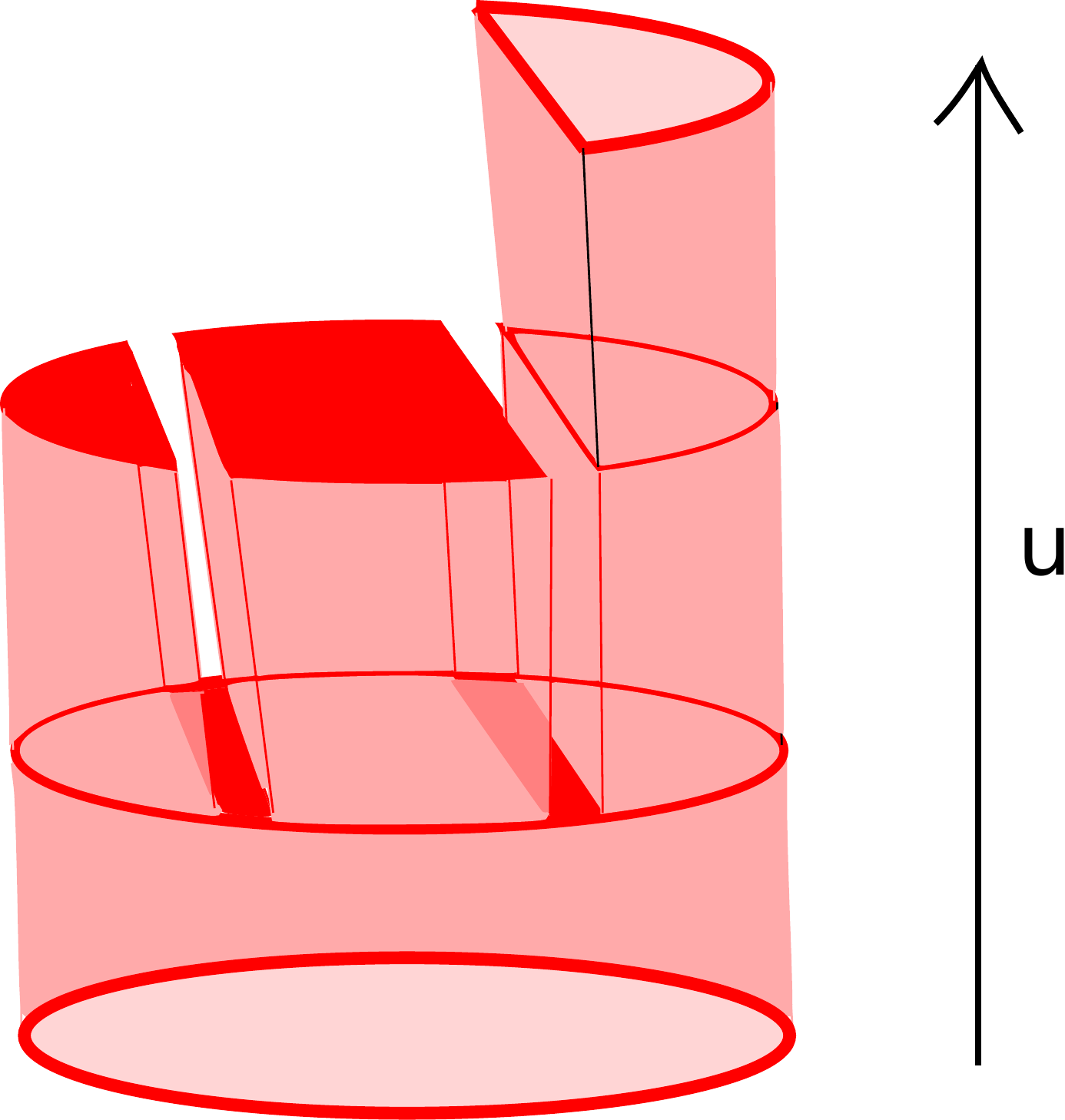}
\caption{This is an illustration of a concordance between circles, but acts as a schematic for the concordance from $R$. We attach $3$-dimensional $1$-handles to a thickened $R$ in $S^2\times S^2\times I$, where the $1$-handles live within one $S^2\times S^2\times u$. We then attach $3$-dimensional $2$-handles to the cobordism (in a later $S^2\times S^2\times u'$), geometrically cancelling the $1$-handles.}
\label{fig:schematic}
\end{figure}

We now perform the construction of Yildiz \cite{yildiz} in the region $-4\le t\le 4$. We will build a cobordism between $R$ and a standard surface in $S^2\times S^2$ by attaching successive handles to (a thickened) $R$ in $S^2\times S^2\times I$, where each handle is contained in one $S^2\times S^2\times u$. The handles will cancel geometrically so that the final cobordism will in fact be a concordance. See Figure \ref{fig:schematic} for a schematic.

\begin{figure}
\includegraphics[width=\textwidth]{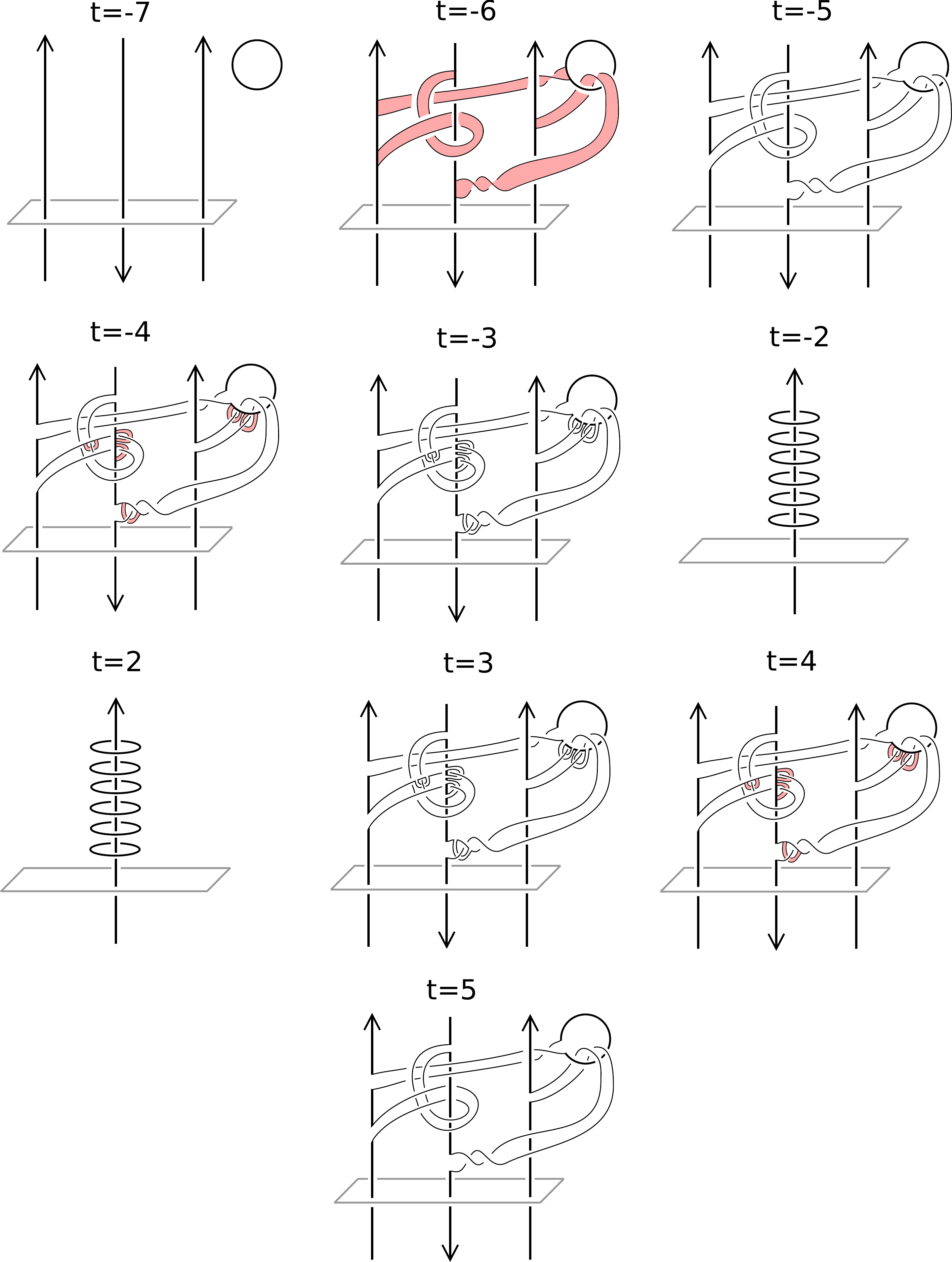}
\caption{We attach $k$ $3$-dimensional $1$-handles to a thickened $R$ to find a cobordism from $R$ to a genus-$(g+k)$ surface $\Sigma$ which intersects $S^2\times S^1\times t$ in the union of $\pt\times S^1$ and $k$ meridians (for $-3<t<3)$. This figure is an illustration of the surface $\Sigma\subset S^2\times S^2$.} 
\label{fig:afteronehandles}
\end{figure}

\begin{figure}
\includegraphics[width=110mm]{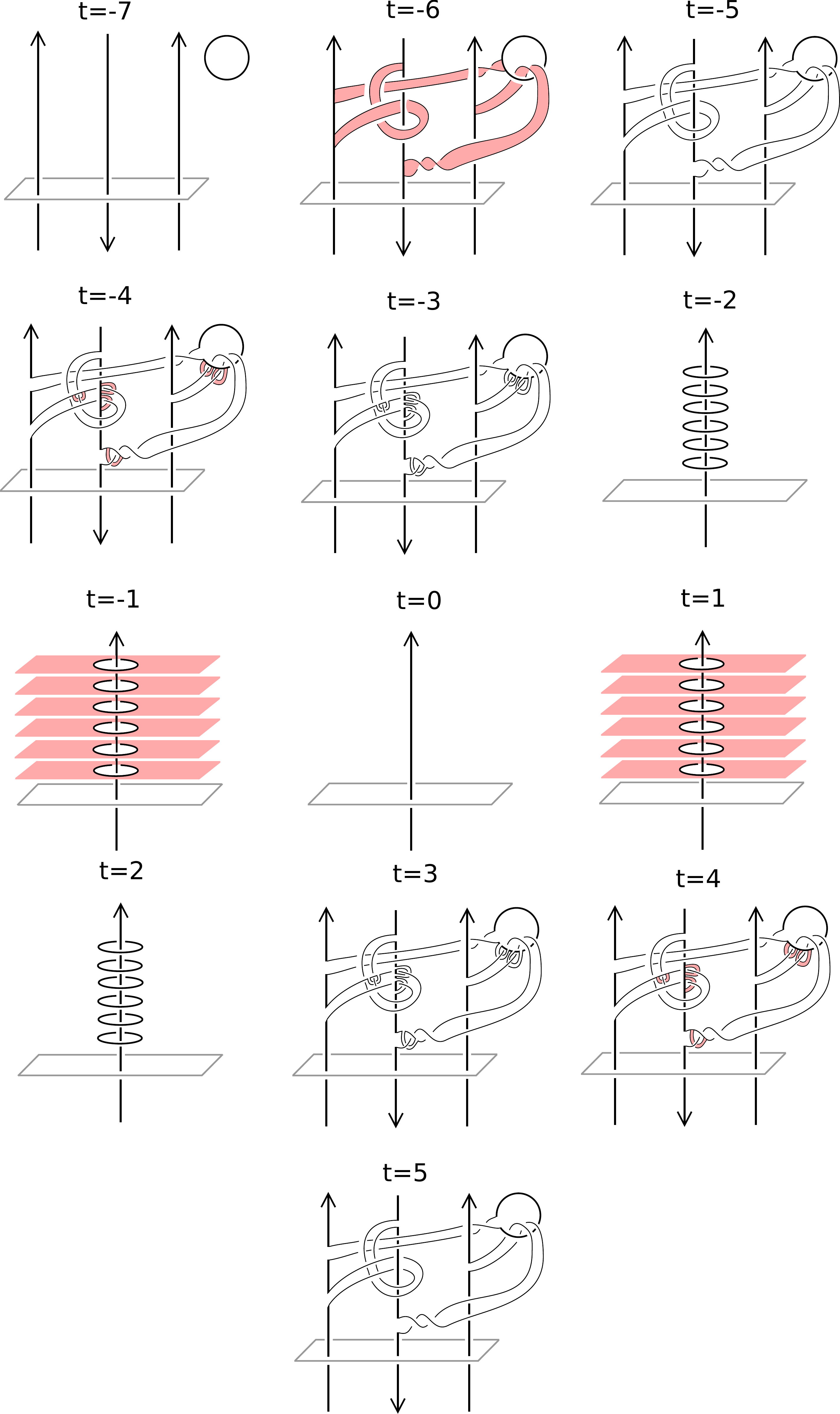}
\caption{We attach $k$ $2$-handles to the cobordism between $R$ and $\Sigma$, geoemetrically cancelling the $1$-handles. This yields a concordance between $R$ and a genus-$g$ surface $R'$ with $R'\cap(S^2\times y_0)=\pt$. This figure is an illustration of the surface $R'$. By \cite{gabai}, $R'$ may be isotoped to be standard.}
\label{fig:afterhandles}
\end{figure}

Attach $k$ $3$-dimensional $1$-handles of the form band$\times[-4,4]$ to a thickened $R$ to find a cobordism from $R$ to a genus-$(g+k)$ surface $\Sigma\subset S^2\times S^2$. Choose the $1$-handles to achieve crossing changes of $R\cap \left(S^2\times S^1\times 0\right)$ (see Fig. \ref{fig:afteronehandles}) so that for $-4<t<4$, $\Sigma\cap(S^2\times S^1\times t)$ is a link of the form $K_t\sqcup_{i=1}^g L_i$, where the $L_i$ are meridians for $K_t$ and $K_t$ is isotopic in $S^2\times S^1\times t$ to a curve meeting $S^2\times (z_0,t)$ once. 


Use the 3-dimensional light bulb theorem to isotope $K_t$ within $-3<t<3$ to be the standard vertical $\pt\times S^1$. See Figure \ref{fig:afteronehandles}.

For $-1\le t\le 1$ (in particular), each $L_i$ bounds a disk $D_i$ in $S^2\times S^1\times t$, where $\int D_i\cap R=\emptyset$. Attach to this cobordism $k$ $3$-dimensional $2$-handles of the form $D_i\times[-1,1]$. 
(See Fig. \ref{fig:afterhandles}.) The $2$-handles geometrically cancel the $1$-handles, yielding a concordance from $R$ to a genus-$g$ surface $R'$ with $R'\cap \left(S^2\times S^1\times 0\right)=K_0$. Therefore, $|R'\cap (S^2\times y_0)=1|$.

Since $R,R'$ are concordant, $[R']=[R]=[R^{\std}]\in H_2(S^2\times S^2;\Z)$. By the 4-dimensional light bulb theorem \cite{gabai}, $R'$ is isotopic to $R^{\std}$ with $g$ standard handles attached. 

\qed

\bibliographystyle{abbrv}

\end{document}